\documentclass[a4paper]{amsart}
\usepackage{amsmath,amssymb,amsfonts}
\usepackage{colordvi}
\usepackage[all]{xy}

\allowdisplaybreaks

\def\mpoint{\;.}
\def\mvirg{\;,}
\def\sur{\overline}
\def\sou{\underline}

\def\mpn{\medskip\par\noindent}

\def\mmpn{\vskip 1em minus 1em\par\noindent}

\def\sp{\bigskip\par}

\def\smp{\smallskip\par}

\def\CC{{\mathcal C}}

\def\CF{{\mathcal F}}

\def\CH{{\mathcal H}}

\def\CL{{\mathcal L}}

\def\CR{{\mathcal R}}

\def\Id{\operatorname{id}\nolimits}

\def\Mod{\operatorname{Mod}\nolimits}
\def\Alg{\operatorname{Alg}\nolimits}
\def\Hom{\operatorname{Hom}\nolimits}

\def\op{^{op}}
\def\dual{^{\scriptscriptstyle\natural}}

\def\meet{\wedge}

\def\N{\mathbb{N}}

\def\S{\mathbb{S}}

\newcommand{\dirsum}[1]{\mathop{\bigoplus}_{#1}\limits}

\newcommand{\edge}[2]{\xymatrix{#1\ar@{->-}[r]&#2}}
\def\marc[#1]{\ar@{-}[#1]|(.4){\object@{<}}}
\def\mard[#1]{\ar@{-}[#1]|(.5){\object@{>}}}
\def\marb[#1]{\ar@{-}[#1]|{\object+{  }}}
\newcommand{\fleche}[2]{\xymatrix@C=4ex{*!U(0.2){#1\;}&*!U(0.5){\;#2}\marc[l]}}
\newcommand{\flecheb}[2]{\xymatrix@C=4ex{*!U(0.2){#1\;}&*!U(0.1){\;#2}\marc[l]}}

\def\pf{\par\bigskip\noindent{\bf Proof~: }}
\def\endpf{~\hfill\rlap{\hspace{-1ex}\raisebox{.5ex}{\framebox[1ex]{}}\sp}\bigskip\pagebreak[3]}

\makeatletter
\renewenvironment{enumerate}{\ifnum \@enumdepth >3 \@toodeep\else
       \advance\@enumdepth \@ne
       \edef\@enumctr{enum\romannumeral\the\@enumdepth}\list
       {\csname  label\@enumctr\endcsname}{\setlength{\topsep}{1ex}
\setlength{\itemsep}{0 pt}\usecounter
         {\@enumctr}\def\makelabel##1{\hss\llap{##1}}}\fi}{\endlist}

\def\@seccntformat#1{\csname the#1\endcsname.\quad}

\def\section{\pagebreak[3]\setcounter{prop}{0}\setcounter{equation}{0}\@startsection{section}{1}{\z@}{4ex plus  6ex}{2ex}{\center\reset@font \large\bf}}

\makeatother

\def\theprop{\thesection.\arabic{prop}}

\renewenvironment{equation}{\refstepcounter{subsection}\refstepcounter
{prop}$$}{\leqno{\bf (\theprop)}$$}

\newenvironment{enonce}[1]{\pagebreak[3]\refstepcounter{prop}\mmpn
{{\bf  \thesection.\arabic{prop}.\ #1.}}\begin{it} }{\end{it}\smp}
\def\thesection{\arabic{section}}

\newcommand{\result}[1]{\begin{enonce}{#1}}
\newcommand{\fresult}{\end{enonce}}

\newcommand{\mbigvee}[1]{\mathop{\bigvee}_{#1}\limits}

\begin{document}

\title[Tensor product of correspondence functors]{Tensor product of correspondence functors}

\author{Serge Bouc}
\author{Jacques Th\'evenaz}



\begin{abstract}
As part of the study of correspondence functors, the present paper investigates their tensor product and proves some of its main properties.
In particular, the correspondence functor associated to a finite lattice has the structure of a commutative algebra in the tensor category of all correspondence functors.
\end{abstract}

\maketitle


\section{Introduction}

\medskip
\noindent
In recent papers \cite{BT2, BT3, BT4}, we developed the theory of correspondence functors, namely functors from the category~$\CC$ of finite sets and correspondences to the category $k\text{-\!}\Mod$ of all $k$-modules, where $k$ is a commutative ring.
This theory turns out to be very rich and we describe here another piece of structure.
We introduce the tensor product of any two correspondence functors (Section~\ref{Section:Tensor}) and analyse its main properties,
such as projectivity, finite generation, and behaviour under induction (Section~\ref{Section:Properties}).\par

Whenever $M$ and~$M'$ are correspondence functors, we not only define their tensor product $M\otimes M'$,
but also their internal hom $\CH(M,M')$, which satisfies the usual adjointness property (Section~\ref{Section:Hom}).
The constructions and its properties depend on the symmetric monoidal structure on $\CC$ given by the disjoint union of finite sets. We also show that our construction turns out to be an instance of the general construction known as Day convolution.\par

A main instance of correspondence functors is the functor $F_T$ associated to a finite lattice~$T$,
as defined in~\cite{BT3}.
We prove in Section~\ref{Section:Lattices} that $F_T\otimes F_{T'}\cong F_{T\times T'}$.
Finally, in Section~\ref{Section:Algebra}, we show that $F_T$ carries the additional structure of a commutative algebra
in the tensor category~$\CF_k$ of all correspondence functors.
We prove conversely that, over a field~$k$ (or more generally if ${\rm Spec}(k)$ is connected),
any commutative algebra in the tensor category~$\CF_k$ is isomorphic to~$F_T$ for some finite lattice~$T$,
provided it satisfies two additional conditions, one of them being of an exponential nature (see Theorem~\ref{converse}).
A few small examples are discussed in Section~\ref{Section:Examples}.\par

Throughout this paper, $k$ denotes a fixed commutative ring and all modules are left $k$-modules.


\bigbreak

\section{Correspondence functors}\label{Section:Basic}

\medskip
\noindent
In this section, we recall the definitions and basic properties of correspondence functors.
We refer to Sections~2--4 of~\cite{BT2} and Section~2 of~\cite{BT3} for more details.
We denote by $\CC$ the category of finite sets and correspondences.
Its objects are the finite sets and the set $\CC(Y,X)$ of morphisms from $X$ to~$Y$ (using a reverse notation which is convenient for left actions)
is the set of all correspondences from $X$ to~$Y$, namely all subsets of~$Y\times X$.
Given two correspondences $V\subseteq Z\times Y$ and $U\subseteq Y\times X$,
their composition $VU$ is defined by
$$VU:=\{ \, (z,x)\in Z\times X \,\mid\, \exists\;y\in Y \;\text{ such that } \; (z,y)\in V \,\text{ and }\, (y,x)\in U \,\} \mpoint$$
The correspondence $\Delta_X=\{(x,x)\mid x\in X\}$ is a unit element for this composition.
A correspondence from $X$ to~$X$ is also called a {\em relation} on~$X$.
In particular,
$$\CR_X:=k\CC(X,X)$$
is a $k$-algebra, the algebra of the monoid of all relations on~$X$.\par

For our fixed commutative ring~$k$, we let $k\CC$ be the $k$-linearization of~$\CC$.
The objects are again the finite sets and $k\CC(Y,X)$ is the free $k$-module with basis $\CC(Y,X)$.
A {\em correspondence functor over~$k$} is a $k$-linear functor from $k\CC$ to the category $k\text{-\!}\Mod$ of all $k$-modules.
We let $\CF_k$ be the category of all correspondence functors over~$k$ (for some fixed commutative ring~$k$).\par

If $M$ is a correspondence functor and $U\in\CC(Y,X)$ is a correspondence, then $U$ acts as a linear map
$M(U):M(X)\to M(Y)$ and we simply write $U$ for this left action.
In other words, for any $m\in M(X)$, we define
$$Um := M(U)(m) \in M(Y) \mpoint$$
In particular, we have $(VU)m=V(Um)$ for any correspondence $V\in\CC(Z,Y)$.\par

The following examples have been introduced in~\cite{BT2} and~\cite{BT3} respectively.

\result{Example} \label{LEW}
{\rm
For any finite set~$E$, the representable functor $k\CC(-,E)$ is a projective correspondence functor by Yoneda's lemma.
More generally, for any left $\CR_E$-module~$W$, there is a correspondence functor $L_{E,W}$
defined by
$$L_{E,W}(X):=k\CC(X,E) \otimes_{\CR_E}  W \mpoint$$
This is left adjoint to the evaluation at~$E$, in the sense that it satisfies the adjointness property
$$\Hom_{\CF_k}(L_{E,W},M) \cong \Hom_{\CR_E}(W,M(E)) \mpoint$$
This example is used in~\cite{BT2},
but it is in fact a general construction for representations of categories which goes back to~\cite{Bo}.
In particular, it is used for the construction of simple functors.
}
\fresult

\result{Example} \label{constant}
{\rm
The constant functor $\sou{k}$ is defined by $\sou{k}(X)=k$ for any finite set~$X$
and $U\lambda=\lambda$ for any $\lambda\in \sou{k}(X)$ and any correspondence $U\in\CC(Y,X)$.
Actually, it is elementary to check that $\sou{k}\cong k\CC(-,\emptyset)$, so in particular $\sou{k}$ is projective.
}
\fresult

\result{Example} \label{F_T}
{\rm For any finite lattice~$T$, let $F_T(X)=kT^X$ for any finite set~$X$,
where $T^X$ is the set of all maps $X\to T$ and $kT^X$ denotes the free $k$-module with basis~$T^X$.
The action of a correspondence $U\in\CC(Y,X)$ on a function $\varphi\in T^X$ is a function $U\varphi \in T^Y$ defined by the join
$$(U\varphi)(y):=\mbigvee{\substack{x\in X \\ (y,x)\in U}} \varphi(x)\mvirg$$
with the usual comment that the join over an empty set yields the least element~$\hat0$ of the lattice.
It is easy to check that this defines a correspondence functor~$F_T$ (see~\cite{BT3}).\par

Recall that a join-morphism $f:T\to T'$ of finite lattices is a map such that $f(\bigvee_{t\in A}t)=\bigvee_{t\in A}f(t)$ for any subset $A$ of~$T$.
In particular, $f$ must be order-preserving.
Moreover, the case when $A$ is empty shows that $f(\hat0)=\hat0$, where $\hat0$ denotes the least element of any lattice.
Any join morphism $f:T\to T'$ induces a morphism of correspondence functors $F_T\to F_{T'}$ by composition with~$f$.
}
\fresult


\bigbreak

\section{Tensor product}\label{Section:Tensor}

\medskip
\noindent
In this section, we define the tensor product of two correspondence functors and discuss its basic properties.

\result{Definition} Let $M$ and $M'$ be correspondence functors over~$k$.
The {\em tensor product} of $M$ and $M'$ is the correspondence functor $M\otimes M'$ defined as follows.
For any finite set~$X$,
$$(M\otimes M')(X)=M(X)\otimes_k M'(X)\mpoint$$
If $Y$ is a finite set and $U\in\CC(Y,X)$, then $U$ acts as a $k$-linear map
$$U:(M\otimes M')(X)\longrightarrow (M\otimes M')(Y)$$ 
defined by
$$U(m\otimes m')=Um\otimes Um' ,\;\;  \forall m\in M(X),\;\forall m'\in M'(X)\mpoint$$
\fresult

The properties of this tensor product use the symmetric monoidal structure on~$\CC$ given by the disjoint union of finite sets,
which we write $\sqcup$  throughout this paper.
More precisely, if $X,X',Y,Y'$ are finite sets, and if $U\in\CC(X',X)$ and $V\in\CC(Y',Y)$,
then the disjoint union $U\sqcup V$ can be viewed as a subset of $(X'\sqcup Y')\times (X\sqcup Y)$.
We will represent this correspondence in the matrix form
$$\left(\begin{array}{cc}U&\emptyset\\\emptyset&V\end{array}\right)\mpoint$$
This yields a bifunctor $\CC\times \CC\to \CC$, still denoted by a disjoint union symbol.\par

We also use the following matrix notation. If $U\in\CC(X',X)$ and $V\in\CC(X'',X)$, then
$$\binom{U}{V}\in \CC(X'\sqcup X'',X)$$
denotes the obvious subset of $(X'\sqcup X'')\times X$.
Similarly, if $U\in\CC(X,X')$ and $V\in\CC(X,X'')$, then
$$(U,V)\in \CC(X,X'\sqcup X'')$$
denotes the obvious subset of $X\times (X'\sqcup X'')$.

\result{Proposition} \label{tensor-product}
Let $M$, $M'$ and $M''$ be correspondence functors over~$k$. 
\begin{enumerate}
\item The assignment $(M,M')\mapsto M\otimes M'$ is a $k$-linear bifunctor $\CF_k\times\CF_k\to \CF_k$ and is right exact in $M$ and~$M'$.
\item There are natural isomorphisms of correspondence functors
\begin{eqnarray*}
M\otimes(M'\otimes M'')&\cong&(M\otimes M')\otimes M''\\
M\otimes M'&\cong&M'\otimes M\\
M\otimes(M'\oplus M'') &\cong& (M\otimes M')\oplus(M\otimes M') \\
\sou{k}\otimes M&\cong& M\mvirg
\end{eqnarray*}
where $\sou{k}$ is the constant functor introduced in Example~\ref{constant}.
\end{enumerate}
\fresult

\pf (a) This is a straightforward consequence of the usual properties of tensor products.\mpn

(b) For any finite set $X$, the standard $k$-linear isomorphisms
\begin{eqnarray*}
M(X)\otimes_k\big(M'(X)\otimes_k M''(X)\big)&\cong&\big(M(X)\otimes_k M'(X)\big)\otimes_k M''(X)\\
M(X)\otimes_k M'(X)&\cong&M'(X)\otimes_k M(X)\\
M(X)\otimes\big(M'(X)\oplus M''(X)\big) &\cong& \big(M(X)\otimes M'(X)\big)\oplus\big(M(X)\otimes M''(X)\big) \\
\sou{k}(X)\otimes_k M(X)=k\otimes_kM(X)&\cong& M(X)\mvirg
\end{eqnarray*}
are clearly compatible with the action of correspondences.
\endpf

There is also a connection between tensor product and bilinear pairings.
If $M$, $M'$, and $M''$ are correspondence functors over~$k$, the $k$-module of {\em bilinear pairings} $M'\times M\to M''$ is the $k$-module of natural transformations from the bifunctor
$$\CC\times\CC \longrightarrow k\hbox{-\rm Mod} \;, \qquad (X,Y)\mapsto M'(X)\otimes_k M(Y)$$
to the bifunctor
$$\CC\times\CC \longrightarrow k\hbox{-\rm Mod} \;, \qquad (X,Y)\mapsto M''(X\sqcup Y) \mpoint$$

\result{Theorem} \label{pairings}
Let $M$, $M'$ and $M''$ be correspondence functors over~$k$. 
The $k$-module of all bilinear pairings $M'\times M\to M''$ is isomorphic to the $k$-module ${\Hom_{\CF_k}(M'\otimes M,M'')}$.
\fresult

\pf
Given a morphism of correspondence functors $\psi:M'\otimes M\to M''$, we need to construct a bilinear pairing $\widehat\psi:M'\times M\to M''$.
For any finite set $X$, there is a $k$-linear map $\psi_X:M'(X)\otimes_kM(X)\to M''(X)$ with the property that,
for any finite set $Z$ and any correspondence $U\in\CC(Z,X)$, the diagram
$$\xymatrix{
M'(X)\otimes_kM(X)\ar[r]^-{\psi_X}\ar[d]_-{U\otimes U}&M''(X)\ar[d]^-{U}\\
M'(Z)\otimes_kM(Z)\ar[r]^-{\psi_Z}&M''(Z)\\
}
$$
is commutative. If $X$ and $Y$ are finite sets, we define a map
$$\widehat{\psi}_{X,Y}:M'(X)\otimes_kM(Y) \longrightarrow M''(X\sqcup Y)$$
as the following composition
$$\xymatrix{
M'(X)\otimes_kM(Y)\ar[rr]^-{\binom{\Delta_X}{\emptyset}\otimes \binom{\emptyset}{\Delta_Y}}&&M'(X\sqcup Y)\otimes_kM(X\sqcup Y)\ar[r]^-{\psi_{X\sqcup Y}}& M''(X\sqcup Y) \mpoint
}
$$
If $X'$ and $Y'$ are finite sets, if $U\in\CC(X',X)$ and $V\in \CC(Y',Y)$, then we claim that the diagram
\newsavebox{\matrice}
\sbox{\matrice}{\begin{footnotesize}$\left(\!\!\begin{array}{cc} U&\!\!\!\!\emptyset\\\emptyset &\!\!\!\!V\end{array} \!\!\right)$\end{footnotesize}}
$$\xymatrix{
M'(X)\otimes_k M(Y)\ar[d]_{U\hspace{4ex}}^{\hspace{4ex}V}\ar[rr]^-{\binom{\Delta_X}{\emptyset}\otimes\binom{\emptyset}{\Delta_Y}}&&M'(X\sqcup Y)\otimes_kM(X\sqcup Y)\ar[d]^{\hspace{2ex}\usebox{\matrice}}_{\usebox{\matrice}\hspace{2ex}}\ar[r]^-{\psi_{X\sqcup Y}}&M''(X\sqcup Y)\ar[d]_{\usebox{\matrice}}\\
M'(X')\otimes_k M(Y')\ar[rr]_-{\binom{\Delta_{X'}}{\emptyset}\otimes\binom{\emptyset}{\Delta_{Y'}}}&&M'(X'\sqcup Y')\otimes_kM(X'\sqcup Y')\ar[r]_-{\psi_{X'\sqcup Y'}}&M''(X'\sqcup Y')\\
}
$$
is commutative.
The left hand side square is commutative because
$$\left(\begin{array}{cc}U&\emptyset\\\emptyset&V\end{array}\right) \binom{\Delta_X}{\emptyset}
=\binom{U}{\emptyset} = \binom{\Delta_{X'}}{\emptyset} U\mvirg$$
and similarly
$$\left(\begin{array}{cc}U&\emptyset\\\emptyset&V\end{array}\right) \binom{\emptyset}{\Delta_Y}
= \binom{\emptyset}{V} = \binom{\emptyset}{\Delta_{Y'}} V \mpoint$$
The right hand side square is commutative by the defining property of the morphism $\psi:M'\otimes M\to M''$.
It follows that
$$\left(\begin{array}{cc}U&\emptyset\\\emptyset&V\end{array}\right)  \widehat{\psi}_{X,Y}
=\widehat{\psi}_{X',Y'} (U\otimes V) \mvirg$$
so that the family of maps $\widehat{\psi}_{X,Y}$ define a natural transformation from the bifunctor $(X,Y)\mapsto M'(X)\otimes_k M(Y)$ to the bifunctor $(X,Y)\mapsto M''(X\sqcup Y)$, in other words a bilinear pairing $\widehat{\psi}:M'\times M\to M''$.\par

Conversely, given a bilinear pairing $\eta:M'\times M\to M''$, we have to construct a morphism of correspondence functors $\widetilde{\eta}:M'\otimes M \to M''$.
For any finite sets $X,Y$, there is a $k$-linear map
$$\eta_{X,Y}:M'(X)\otimes_kM(Y) \longrightarrow M''(X\sqcup Y)$$
such that, for any finite set $X'$ and any correspondences $U\in\CC(X',X)$ and $V\in\CC(Y',Y)$, the diagram
\begin{equation}\label{diagram}
\vcenter{\xymatrix{
M'(X)\otimes_kM(Y)\ar[d]_{U\otimes V}\ar[r]^-{\eta_{X,Y}}&M''(X\sqcup Y)\ar[d]^{\usebox{\matrice}}\\
M'(X')\otimes_kM(Y')\ar[r]^-{\eta_{X',Y'}}&M''(X'\sqcup Y')\\
}}
\end{equation}
is commutative.\par

In particular, for $X=Y$, we have a map $\eta_{X,X}:M'(X)\otimes_kM(X)\to M''(X\sqcup X)$ which we can compose with the map $M''(X\sqcup X)\to M''(X)$ given by the action of the correspondence $(\Delta_X,\Delta_X)\in\CC(X,X\sqcup X)$, to get a map
$$\widetilde{\eta}_{X} := (\Delta_X,\Delta_X) \, \eta_{X,X}: M'(X)\otimes_k M(X) \longrightarrow M''(X)\mpoint$$
\sbox{\matrice}{\begin{footnotesize}$\left(\!\!\begin{array}{cc}U &\!\!\!\!\emptyset\\ \emptyset&\!\!\!\! U\end{array}\!\!\right)$\end{footnotesize}}
If $Z$ is a finite set and $U\in\CC(Z,X)$, we claim that the diagram
$$\xymatrix{
M'(X)\otimes_k M(X)\ar[d]_{U\otimes U}\ar[r]^-{\eta_{X,X}}
&M''(X\sqcup X)\;\;\ar[d]^{\usebox{\matrice}}\ar[r]^-{(\Delta_X,\Delta_X)}& \;M''(X)\ar[d]^U\\
M'(Z)\otimes_k M(Z)\ar[r]_-{\eta_{Z,Z}}&M''(Z\sqcup Z)\ar[r]_-{(\Delta_Z,\Delta_Z)}&M''(Z)\\
}
$$
is commutative. The commutativity of the left hand side square is a special case of the commutativity of the diagram~(\ref{diagram}).
The right hand side square is commutative because
$$U(\Delta_X,\Delta_X)=(U,U)=(\Delta_Z,\Delta_Z)\usebox{\matrice} \mpoint$$
It follows that
$$U \widetilde{\eta}_{X} = \widetilde{\eta}_{Z} (U\otimes U)$$
and therefore the family of maps $\widetilde{\eta}_X$ define a morphism of correspondence functors $\widetilde{\eta}$ from $M'\otimes M$ to $M''$.\par

The constructions $\psi\mapsto\widehat{\psi}$ and $\eta\mapsto\widetilde{\eta}$ are $k$-linear and we need to prove that they are inverse to each other.\par

Let $\psi:M'\otimes M\to M''$ be a morphism of correspondence functors.
For any finite set $X$, we have
\begin{eqnarray*}
\widetilde{\widehat{\psi}\,}_X= (\Delta_X,\Delta_X) \, \widehat{\psi}_{X,X}
&=& (\Delta_X,\Delta_X) \, \psi_{X\sqcup X} \Big( \binom{\Delta_X}{\emptyset} \otimes \binom{\emptyset}{\Delta_X} \Big)\\
&=&\psi_X \Big( (\Delta_X,\Delta_X) \otimes (\Delta_X,\Delta_X) \Big) \Big( \binom{\Delta_X}{\emptyset} \otimes \binom{\emptyset}{\Delta_X} \Big)\\
&=&\psi_X (\Delta_X\otimes\Delta_X) \mvirg
\end{eqnarray*}
because $(\Delta_X,\Delta_X) \binom{\Delta_X}{\emptyset}=\Delta_X$
and $(\Delta_X,\Delta_X) \binom{\emptyset}{\Delta_X}=\Delta_X$.
Since $\Delta_X\otimes\Delta_X$ acts as the identity on $M'(X)\otimes_k M(X)$, we get $\widetilde{\widehat{\psi}\,}_X=\psi_X$, as required.\par

Now let $\eta:M'\times M\to M''$ be a bilinear pairing.
For any finite sets $X$ and $Y$, the definition of $\widehat \psi_{X,Y}$ applied to $\psi=\widetilde\eta$ and the definition of~$\widetilde\eta_{X\sqcup Y}$ yield
\begin{eqnarray*}
\widehat{\;\widetilde{\eta}\;}_{X,Y}&=&\widetilde{\eta}_{X\sqcup Y}
\Big( \binom{\Delta_X}{\emptyset} \otimes \binom{\emptyset}{\Delta_Y} \Big) \\
&=& (\Delta_{X\sqcup Y},\Delta_{X\sqcup Y}) \, \eta_{X\sqcup Y,X\sqcup Y} \Big( \binom{\Delta_X}{\emptyset} \otimes \binom{\emptyset}{\Delta_Y}\Big) \\
&=& (\Delta_{X\sqcup Y},\Delta_{X\sqcup Y})
\left({\begin{array}{cc}\scriptstyle\binom{\Delta_X}{\emptyset}&\scriptstyle\emptyset \\
\scriptstyle\emptyset&\scriptstyle\binom{\emptyset}{\Delta_Y}\end{array}}\right) \eta_{X,Y} \mvirg
\end{eqnarray*}
the latter equality coming from the commutative diagram~(\ref{diagram}) for the sets $X'=Y'=X\sqcup Y$
and the correspondences $U=\binom{\Delta_X}{\emptyset}$ and $V=\binom{\emptyset}{\Delta_Y}$.
Now it is easy to check that
$$(\Delta_{X\sqcup Y},\Delta_{X\sqcup Y})
\left({\begin{array}{cc}\scriptstyle\binom{\Delta_X}{\emptyset}&\scriptstyle\emptyset \\
\scriptstyle\emptyset&\scriptstyle\binom{\emptyset}{\Delta_Y}\end{array}}\right)
=\left({\begin{array}{cc} \Delta_X&\emptyset \\ \emptyset&\Delta_Y\end{array}}\right) \mvirg$$
and this acts as the identity on~$M'(X)\otimes_k M(Y)$.
Therefore $\widehat{\;\widetilde{\eta}\;}_{X,Y}=\eta_{X,Y}$, as was to be shown.
\endpf

Theorem~\ref{pairings} shows that our construction of tensor product is a special case of the general construction due to Day, known as Day convolution  \cite{Da}.


\bigbreak

\section{Lattices}\label{Section:Lattices}

\medskip
\noindent
We want to apply the tensor product construction to functors of the form~$F_T$, where $T$ is a finite lattice, as defined in Example~\ref{F_T}.
As in~\cite{BT3}, we define the category~$k\CL$ of finite lattices as follows.
Its objects are the finite lattices and $\Hom_{k\CL}(T,T')$ is the free $k$-module with basis the set of all join-morphisms from $T$ to~$T'$.\par

The direct product $T\times T'$ of two lattices is defined using componentwise operations.
Our next lemma shows that there is also a direct product for morphisms in~$k\CL$.

\result{Lemma} \label{product-morphisms} Let $S$, $T$, S'$, T'$ be finite lattices.
\begin{enumerate}
\item If $f:S\to T$ and $f':S'\to T'$ are join-morphisms, then $f\times f': S\times S' \to T\times T'$ is a join-morphism.
\item Extending this direct product by $k$-bilinearity defines a $k$-linear bifunctor $k\CL \times k\CL \to k\CL$.
\end{enumerate}
\fresult

\pf
Given a subset $A\subseteq S\times S'$, let $B\subseteq S$ (respectively $B'\subseteq S'$) denote the first projection of $A$ (respectively the second projection).
Then
\begin{eqnarray*}
(f\times f')\big(\mbigvee{(s,s')\in A}(s,s')\big)&=&(f\times f')\Big(\mbigvee{(s,s')\in A}\big((s,\hat0)\vee(\hat0,s')\big)\Big)\\
&=&(f\times f')\Big(\big(\mbigvee{s\in B}(s,\hat0)\big) \vee \big(\mbigvee{s'\in B'}(\hat0,s')\big) \Big)\\
&=&(f\times f')\big(\mbigvee{s\in B}s\,,\mbigvee{s'\in B'}s'\big) \\
&=&\Big(f(\mbigvee{s\in B} s)\,,f'(\mbigvee{s'\in B'}s')\Big)\\
&=&\Big(\mbigvee{s\in B}f(s),\mbigvee{s'\in B'}f'(s')\Big)\\
&=&\Big(\mbigvee{s\in B}f(s)\,,\hat0\Big)\vee\Big(\hat0,\mbigvee{s'\in B'}f'(s')\Big)\\
&=&\Big(\mbigvee{s\in B}\big(f(s),\hat{0}\big) \Big) \vee \Big(\mbigvee{s'\in B'}\big(\hat0,f'(s')\big) \Big)\\
&=&\mbigvee{(s,s')\in A}\big(f(s),f'(s')\big)
\end{eqnarray*}
as required. The second assertion follows by bilinearity.
\endpf

\result{Theorem} \label{tensor lattice} 
The bifunctors
$$k\CL\times k\CL \longrightarrow \CF_k \,,\qquad (T,T')\mapsto F_T\otimes F_{T'}$$
and
$$k\CL\times k\CL \longrightarrow \CF_k \,,\qquad (T,T')\mapsto F_{T\times T'}$$ 
are isomorphic.
\fresult

\pf 
Let $T$ and $T'$ be finite lattices. For any finite set $X$, there is a unique isomorphism of $k$-modules
$$\tau_X:(F_T\otimes F_{T'})(X)=k(T^X)\otimes_k k(T'^X)\longrightarrow k(T\times T')^X\mvirg$$
mapping $\varphi\otimes\varphi'$ to $\varphi\times\varphi'$.
Here of course, the map $\varphi\times\varphi':X\to T\times T'$ is obtained by direct product form the maps $\varphi:X\to T$ and $\varphi':X\to T'$.
If $Y$ is a finite set and $U\in \CC(Y,X)$, then for any $y\in Y$,
\begin{eqnarray*}
U(\varphi\times\varphi')(y)&=&\mbigvee{(y,x)\in U}\big(\varphi(x),\varphi'(x)\big)\\
&=&\big(\mbigvee{(y,x)\in U}\varphi(x),\mbigvee{(y,x')\in U}\varphi'(x')\big)\\
&=&\big(U\varphi(y),U\varphi'(y)\big)\mpoint
\end{eqnarray*}
Thus $U(\varphi\times\varphi')=U\varphi\times U\varphi'$, that is,
$U\,\tau_X(\varphi\otimes \varphi')=\tau_Y(U\varphi\otimes U\varphi')$.
Therefore
$$\tau:F_T\otimes F_{T'}\longrightarrow F_{T\times T'}$$
is an isomorphism of correspondence functors.\par

If $f:S\to T$ and $f':S'\to T'$ are join-morphisms, then $f\times f':S\times S'\to T\times T'$ is a  join-morphism, by Lemma~\ref{product-morphisms}.
Moreover, we claim that the diagram
$$\xymatrix{
F_S\otimes F_{S'}\ar[d]_{F_f\otimes F_{f'}}\ar[r]^{\sigma}&F_{S\times S'}\ar[d]^{F_{f\times f'}}\\
F_T\otimes F_{T'}\ar[r]^{\tau}&F_{T\times T'}\\
}
$$
is commutative, where $\sigma:F_S\otimes F_{S'}\to F_{S\times S'}$ denotes the corresponding isomorphism for the lattices $S$ and $S'$.
This is because, for any finite set $X$, any map $\varphi:X\to S$, and any map $\varphi':X'\to S'$, we have
$$F_{f\times f'}\,\sigma_X(\varphi\otimes\varphi')=F_{f\times f'}(\varphi\times\varphi')=(f \varphi)\times (f' \varphi')
=\tau_X \big(F_f(\varphi)\otimes F_{f'}(\varphi')\big)\mpoint$$
It follows that the family of isomorphisms $\tau$ yields an isomorphism between the two bifunctors of the statement.
\endpf

\result{Corollary} \label{tensor representable}
If $E$ and $E'$ are finite sets, then 
$$k\CC(-,E)\otimes k\CC(-,E')\cong k\CC(-,E\sqcup E')\mpoint$$
\fresult

\pf This follows from Theorem~\ref{tensor lattice} applied to the lattice $T$ of subsets of~$E$ and the lattice $T'$ of subsets of~$E'$.
Then $F_T\cong k\CC(-,E)$ and $F_{T'}\cong k\CC(-,E')$, because a map from $X$ to~$T$ corresponds to a subset of $X\times E$.
Moreover $T\times T'$ is isomorphic to the lattice of subsets of~$E\sqcup E'$.
\endpf


\bigbreak

\section{More properties of tensor product}\label{Section:Properties}

\medskip
\noindent
We first discuss projectivity.
Recall that any correspondence functor $M$ is isomorphic to a quotient of $\dirsum{i\in I} k\CC(-,E_i)$ where each $E_i$ is a finite set and $I$ is some index set. This is because if $m_i\in M(E_i)$, Yoneda's lemma gives a morphism $\psi_i: k\CC(-,E_i) \to M$ mapping $\Delta_{E_i}$ to~$m_i$.
Choosing a set $\{m_i\mid i\in I\}$ of generators of~$M$, the sum of the morphisms~$\psi_i$ yields a surjective morphism
$\bigoplus_{i\in I} k\CC(-,E_i) \to M$, as required.
In particular, any projective correspondence functor is isomorphic to a direct summand of a direct sum of representable functors.

\result{Proposition} \label{projective-tensor} Let $M$ and $N$ be correspondence functors over~$k$.
If $M$ and $N$ are projective, then so is $M\otimes N$.
\fresult

\pf
By the observation above, it suffices to assume that $M=\dirsum{i\in I}k\CC(-,E_i)$ and $N=\dirsum{j\in J}k\CC(-,F_j)$, where
$E_i$ and $F_j$ are finite sets and where $I$ and $J$ are some index sets.
By Corollary~\ref{tensor representable}, we obtain
$$M\otimes N\cong\dirsum{i\in I,\,j\in J}k\CC(-,E_i)\otimes k\CC(-,F_j)\cong\dirsum{i\in I,\,j\in J}k\CC(-,E_i\sqcup F_j)\mvirg$$
so $M\otimes N$ is projective.
\endpf

It should be observed that, since $\sou{k}\otimes M\cong M$ for any correspondence functor~$M$ (Proposition~\ref{tensor-product})
and since $\sou{k}$ is projective (Example~\ref{constant}), 
tensoring with a projective functor does not yield a projective functor in general, contrary to the case of finite group representations.\par

Next we consider the functors $L_{E,V}$ introduced in Example~\ref{LEW}, where $E$ is a finite set and $V$ is an $\CR_E$-module.
Recall that $\CR_E:=k\CC(E,E)$.
There is an induction procedure
$$V\!\uparrow_E^F \;:=k\CC(F,E)\otimes_{\CR_E}V \mvirg$$
where $F$ is any finite set. Clearly $V\!\uparrow_E^F$ is a left $\CR_F$-module.
Notice that we have $L_{E,V}(F)=V\!\uparrow_E^F$ by the definition of~$L_{E,V}$.

\result{Theorem} \label{LEV-tensor} Let $E$ and $F$ be finite sets and let $G=E\sqcup F$.
Let $V$ be a $\CR_E$-module and $W$ be a $\CR_F$-module.
Then there is an isomorphism of correspondence functors
$$L_{E,V}\otimes L_{F,W}\cong L_{G,V\uparrow_E^G\otimes_kW\uparrow_F^G}\mvirg$$
where the $k\CC(G,G)$-module structure on $V\uparrow_E^G\otimes_kW\uparrow_F^G$ is induced by the comultiplication
$\nu:k\CC(G,G)\to k\CC(G,G)\otimes_k k\CC(G,G)$ defined by $\nu(U)=U\otimes U$ for any relation $U\in\CC(G,G)$.
\fresult

\pf
Since $(L_{E,V}\otimes L_{F,W})(G)=L_{E,V}(G)\otimes_k L_{F,W}(G)={V\!\uparrow_E^G}\otimes_k{W\!\uparrow_F^G}$,
the identity map ${V\!\uparrow_E^G}\otimes_k{W\!\uparrow_F^G} \to (L_{E,V}\otimes L_{F,W})(G)$ corresponds,
by the adjointness property of $L_{G,V\!\uparrow_E^G\otimes_kW\!\uparrow_F^G}$, to a morphism
$$\Phi: L_{G,V\!\uparrow_E^G\otimes_kW\!\uparrow_F^G} \longrightarrow L_{E,V}\otimes L_{F,W}$$
which we need to described explicitly.
For a finite set $X$, 
$$L_{G,V\uparrow_E^G\otimes_kW\uparrow_F^G}(X)=k\CC(X,G)\otimes_{\CR_G}\Big(\big(k\CC(G,E)\otimes_{\CR_E}V\big)\otimes_k\big(k\CC(G,F)\otimes_{\CR_F}W\big)\Big)$$
and 
$$
(L_{E,V}\otimes L_{F,W})(X)=\big(k\CC(X,E)\otimes_{\CR_E}V\big)\otimes_k\big(k\CC(X,F)\otimes_{\CR_F}W\big)\mpoint$$
It is easy to check that the morphism $\Phi_X$ maps the element 
$$C\otimes_{\CR_G}\big((A\otimes_{\CR_E} v)\otimes_k(B\otimes_{\CR_F} w)\big) \;
\in L_{G,V\uparrow_E^G\otimes_kW\uparrow_F^G}(X)$$
to the element 
$$(CA\otimes_{\CR_E} v)\otimes_k (CB\otimes_{\CR_F} w) \;
\in (L_{E,V}\otimes L_{F,W})(X) \mvirg$$
where $C\in\CC(X,G)$, $A\in\CC(G,E)$, $v\in V$, $B\in\CC(G,F)$, and $w\in W$.\par

Conversely, there is a morphism
$$\Psi_X:(L_{E,V}\otimes L_{F,W})(X)\to L_{G,V\uparrow_E^G\otimes_kW\uparrow_F^G}(X)$$
defined as follows.
For any $P\in\CC(X,E)$, $v\in V$, $Q\in\CC(X,F)$, and $w\in W$, it maps the element
$$(P\otimes_{\CR_E} v)\otimes_k (Q\otimes_{\CR_F} w) \;\in (L_{E,V}\otimes L_{F,W})(X)$$
to the element
$$(P,Q)\otimes_{\CR_G} \Big(\big(\binom{\Delta_E}{\emptyset}\otimes_{\CR_E} v\big)\otimes_k\big(\binom{\emptyset}{\Delta_F}\otimes_{\CR_F} w\big)\Big)
\;\in L_{G,V\uparrow_E^G\otimes_kW\uparrow_F^G}(X) \mvirg$$
where $(P,Q)\in\CC(X,E\sqcup F)$, and where
$\binom{\Delta_E}{\emptyset}\in\CC(E\sqcup F,E)$ and $\binom{\emptyset}{\Delta_F}\in\CC(E\sqcup F,F)$.\par

The map $\Psi_X$ is well defined, for if $R\in\CR_E$ and $S\in\CR_F$, then
\begin{eqnarray*}
&&\Psi_X\big((P\otimes_{\CR_E} Rv)\otimes_k (Q\otimes_{\CR_F} Sw)\big) = \\
&=&(P,Q)\otimes_{\CR_G}
\Big(\big(\binom{\Delta_E}{\emptyset}\otimes_{\CR_E} Rv\big)\otimes_k\big(\binom{\emptyset}{\Delta_F}\otimes_{\CR_F} Sw\big)\Big) \\
&=&(P,Q)\otimes_{\CR_G}
\Big(\big(\binom{R}{\emptyset}\otimes_{\CR_E} v\big)\otimes_k\big(\binom{\emptyset}{S}\otimes_{\CR_F} w\big)\Big) \\
&=&(P,Q)\otimes_{\CR_G}
\left(\begin{array}{cc}R&\emptyset\\\emptyset&S\end{array}\right)
\Big(\big(\binom{\Delta_E}{\emptyset}\otimes_{\CR_E} v\big)\otimes_k\big(\binom{\emptyset}{\Delta_F}\otimes_{\CR_F} w\big)\Big) \\
&=&(PR,QS)\Big(\big(\binom{\Delta_E}{\emptyset}\otimes_{\CR_E} v\big)\otimes_k\big(\binom{\emptyset}{\Delta_F}\otimes_{\CR_F} w\big)\Big) \\
&=&\Psi_X\big((PR\otimes_{\CR_E} v)\otimes_k(QS\otimes_{\CR_F} w)\big) \mpoint
\end{eqnarray*}
Moreover, if $Y$ is a finite set and $U\in\CC(Y,X)$, then
\begin{eqnarray*}
&&\Psi_X\Big(U\big((P\otimes_{\CR_E} Rv)\otimes_k(Q\otimes_{\CR_F} Sw)\big)\Big) =  \\
&=& \Psi_X\big((UP\otimes_{\CR_E} Rv)\otimes_k(UQ\otimes_{\CR_F} Sw)\big) \\
&=& (UP,UQ)\otimes_{\CR_G}
\Big(\big(\binom{\Delta_E}{\emptyset}\otimes_{\CR_E} v\big)\otimes_k\big(\binom{\emptyset}{\Delta_F}\otimes_{\CR_F} w\big)\Big) \\
&=& U\Psi_X\big((P\otimes_{\CR_E} Rv)\otimes_k(Q\otimes_{\CR_F} Sw)\big) \mpoint
\end{eqnarray*}
It follows that the maps $\Psi_X$ define a morphism of correspondence functors 
$$\Psi:L_{E,V}\otimes L_{F,W}\to L_{G,V\uparrow_E^G\otimes_kW\uparrow_F^G}\mpoint$$
Moreover, setting $u=(P\otimes_{\CR_E} v)\otimes(Q\otimes_{\CR_F} w)$, we have
\begin{eqnarray*}
\Phi_X\Psi_X(u)&=&\Phi_X\left((P,Q)\otimes_{\CR_G} \Big(\big(\binom{\Delta_E}{\emptyset}\otimes_{\CR_E} v\big)\otimes_k\big(\binom{\emptyset}{\Delta_F}\otimes_{\CR_F} w\big)\Big)\right)\\
&=&\Big((P,Q)\binom{\Delta_E}{\emptyset}\otimes_{\CR_E}v\Big)\otimes_k\Big((P,Q)\binom{\emptyset}{\Delta_F}\otimes_{\CR_F}w\Big)\\
&=&(P\otimes_{\CR_E} v)\otimes(Q\otimes_{\CR_F} w)=u\mvirg
\end{eqnarray*}
so $\Phi\Psi$ is equal to the identity morphism.\par

Similarly, setting $s=C\otimes_{\CR_G}\big((A\otimes_{\CR_E} v)\otimes_k(B\otimes_{\CR_F} w)\big)$,
\begin{eqnarray*}
\Psi_X\Phi_X(s)&=&\Psi_X\big((CA\otimes_{\CR_E} v)\otimes_k (CB\otimes_{\CR_F} w)\big)\\
&=&(CA,CB)\otimes_{\CR_G} \Big(\big(\binom{\Delta_E}{\emptyset}\otimes_{\CR_E} v\big)\otimes_k\big(\binom{\emptyset}{\Delta_F}\otimes_{\CR_F} w\big)\Big)\\
&=&C(A,B)\otimes_{\CR_G} \Big(\big(\binom{\Delta_E}{\emptyset}\otimes_{\CR_E} v\big)\otimes_k\big(\binom{\emptyset}{\Delta_F}\otimes_{\CR_F} w\big)\Big)\\
&=&C\otimes_{\CR_G}(A,B)\Big(\big(\binom{\Delta_E}{\emptyset}\otimes_{\CR_E} v\big)\otimes_k\big(\binom{\emptyset}{\Delta_F}\otimes_{\CR_F} w\big)\Big)\\
&=&C\otimes_{\CR_G}\big((A\otimes_{\CR_E} v)\otimes_k(B\otimes_{\CR_F} w)\big)=s\mvirg\\
\end{eqnarray*}
so $\Psi\Phi$ is also equal to the identity morphism.
\endpf

Finally we consider finite generation.
Recall from~\cite{BT2} that a correspondence functor~$M$ has {\em bounded type} if there is a finite set~$E$ such that $M$ is generated by~$M(E)$, that is,
$M(X)=k\CC(X,E)M(E)$ for every finite set~$X$.
Moreover, $M$ is {\em finitely generated} if there is a finite set~$E$ and a finite subset $A$ of~$M(E)$ such that $M$ is generated by~$A$, that is,
$M(X)=k\CC(X,E)A$ for every finite set~$X$ (see Proposition~6.4 in~\cite{BT2}).

\result{Theorem} \label{finite-gen} Let $M$ and $N$ be correspondence functors over~$k$.
\begin{enumerate} 
\item If $M$ and $N$ have bounded type, so has $M\otimes N$.
\item If $M$ and $N$ are finitely generated, so is $M\otimes N$.
\end{enumerate}
\fresult

\pf (a) 
Let $E$ and $F$ be finite sets such that $M$ is generated by~$M(E)$ and $N$ is generated by~$N(F)$.
Then the counit morphisms $L_{E,M(E)}\to M$ and $L_{F,N(F)}\to N$ are surjective.
Therefore $M\otimes N$ is isomorphic to a quotient of $L_{E,M(E)}\otimes L_{F,N(F)}$.
By Theorem~\ref{LEV-tensor},
$$L_{E,M(E)}\otimes L_{F,N(F)}\cong L_{G,M(E)\uparrow_E^G\otimes_k N(F)\uparrow_F^G}$$
where $G=E\sqcup F$.
Since $L_{G,M(E)\uparrow_E^G\otimes_k N(F)\uparrow_F^G}$ is generated by its evaluation at~$G$, so is $M\otimes N$.\mpn

(b) Assume now that $M$ is generated by a finite subset $A$ of~$M(E)$ and $N$ is generated by a finite subset $B$ of~$N(F)$.
In particular, $M(E)$ is generated by~$A$ as an $\CR_E$-module and $N(F)$ is generated by~$B$ as an $\CR_F$-module.
Therefore, as $\CR_G$-modules, $M(E)\uparrow_E^G$ is generated by the finite set~$\CC(G,E)\otimes_{\CR_E}A$
and $N(F)\uparrow_F^G$ is generated by the finite set~$\CC(G,F)\otimes_{\CR_F}B$, where $G=E\sqcup F$ as before.
It follows that $M(E)\uparrow_E^G\otimes_k N(F)\uparrow_F^G$ is generated as an $\CR_G$-module by the finite set
$$S:=\big( \CC(G,E)\otimes_{\CR_E}A \big) \otimes_k \big(\CC(G,F)\otimes_{\CR_F}B \big) \mpoint$$
Since $L_{G,M(E)\uparrow_E^G\otimes_k N(F)\uparrow_F^G}$ is generated by its evaluation at~$G$,
namely the $k$-module $M(E)\uparrow_E^G\otimes_k N(F)\uparrow_F^G$, it is also generated by the finite set~$S$.
Now $M\otimes N$ is isomorphic to a quotient of $L_{G,M(E)\uparrow_E^G\otimes_k N(F)\uparrow_F^G}$, so it is generated by the image of~$S$.
Thus $M\otimes N$ is finitely generated.
\endpf


\bigbreak

\section{Internal hom}\label{Section:Hom}

\medskip
\noindent
In this section, we use the symmetric monoidal structure of the category~$\CC$ to define an internal hom in the category~$\CF_k$ of correspondence functors.
We first introduce a useful construction.

\result{Definition} Let $E$ be a finite set and let $M$ be a correspondence functor over~$k$.
\begin{enumerate} 
\item We let $\, t_E: k\CC\to k\CC$ be the endofunctor defined on objects by $t_E(X)=X\sqcup E$ and on correspondences $U\in\CC(Y,X)$ by
$$t_E(U)= U\sqcup \Id = \left(\begin{array}{cc} U&\emptyset\\ \emptyset&\Delta_E \end{array}\right)  \mpoint$$
\item We denote by~$M_E$ the correspondence functor obtained from $M$ by precomposition with the endofunctor $t_E: k\CC\to k\CC$.
\item Let $F$ be a finite set and $V\in\CC(F,E)$. We define $M_V:M_E\to M_F$ to be the morphism obtained by precomposition with the natural transformation $\Id\sqcup V:t_E\to t_F$.
\end{enumerate}
\fresult

Explicitly, we see that $M_E(X)=M(X\sqcup E)$ and $M_V:M(X\sqcup E) \to M(X\sqcup F)$ is given by the action of the correspondence
$$\Id\sqcup V = \left(\begin{array}{cc}\Delta_X&\emptyset\\\emptyset&V\end{array}\right)\mpoint$$

\result{Definition} Let $M$ and $M'$ be correspondence functors over~$k$.
We denote by~$\CH(M,M')$ the correspondence functor defined on a finite set~$E$ by
$$\CH(M,M')(E)=\Hom_{\CF_k}(M,M'_E)\mvirg$$
and for $V\in\CC(F,E)$, by composition with $M'_V:M'_E\to M'_F$.
\fresult

\result{Lemma} \label{internal-hom}
The assignment $(M,M')\mapsto\CH(M,M')$ is a  $k$-linear bifunctor $\CF_k\op\times \CF_k\to \CF_k$, left exact in $M$ and~$M'$.
\fresult

\pf
This is straightforward.
\endpf

Now we prove the basic adjointness property which shows that $\CH(M,M')$ is an internal hom in the category~$\CF_k$.

\result{Theorem} \label{adjointness}
There are isomorphisms of $k$-modules
$$\Hom_{\CF_k}(M'\otimes M,M'')\cong\Hom_{\CF_k}\big(M,\CH(M',M'')\big)$$
natural in $M,M',M''$.
In particular, for any correspondence functor $M'$ over~$k$, the endofunctor
$$\CF_k\longrightarrow \CF_k \,, \qquad M\mapsto M'\otimes M$$
is left adjoint to the endofunctor
$$\CF_k\longrightarrow \CF_k \,, \qquad M\mapsto \CH(M',M) \mpoint$$
\fresult

\pf
Let $\psi: M'\otimes M\to M''$ be a morphism of correspondence functors.
By Theorem~\ref{pairings}, we get a bilinear pairing $M'\times M\to M''$,
hence, for any finite sets $X$ and $Y$, a $k$-linear map
$$\widehat{\psi}_{X,Y}:M'(X)\otimes_k M(Y)\to M''(X\sqcup Y) \mvirg$$
or equivalently, a $k$-linear map 
$$\sur{\psi}_{Y,X}: M(Y)\to \Hom_k\big(M'(X),M''(X\sqcup Y)\big)$$
defined by $\sur{\psi}_{Y,X}(m)(m')=\widehat{\psi}_{X,Y}(m'\otimes m)$, for $m\in M(Y)$ and $m'\in M'(X)$.\par

Now $M''(X\sqcup Y)=M''_Y(X)$.
Moreover, for any finite set $X'$ and any $U\in\CC(X',X)$, the commutative diagram~(\ref{diagram}), for $Y'=Y$ and $V=\Delta_Y$, becomes 
\sbox{\matrice}{\begin{footnotesize}$\left(\!\!\begin{array}{cc}U&\!\!\!\!\emptyset\\ \emptyset&\!\!\!\!\Delta_Y\end{array}\!\!\!\right)$\end{footnotesize}}
$$\xymatrix{
M'(X)\otimes_kM(Y)\ar[d]_{U\otimes\Delta_Y}\ar[r]^-{\widehat{\psi}_{X,Y}}
&M''(X\sqcup Y)\quad \ar[d]^{\usebox{\matrice}} \ar[r]^-{=} & \; M''_Y(X) \ar[d]^{U} \\
M'(X')\otimes_kM(Y)\ar[r]^-{\widehat{\psi}_{X',Y}} &M''(X'\sqcup Y)\quad \ar[r]^-{=} & \; M''_Y(X') \\
}$$
or in other words $\sur{\psi}_{Y,X'}(m)(Um')=U\sur{\psi}_{Y,X}(m)(m')$ for any $m\in M(Y)$ and $m'\in M'(X)$.
Therefore, for a fixed set~$Y$ and a fixed $m\in M(Y)$, the maps $\sur{\psi}_{Y,X}(m)$ define a morphism of correspondence functors
$$\sur{\psi}_Y(m):M'\longrightarrow M''_Y \mvirg$$
hence an element of~$\CH(M',M'')(Y)$.
Allowing $m$ to vary, we obtain a $k$-linear map
$$\sur{\psi}_Y:M(Y)\longrightarrow \CH(M',M'')(Y) \mpoint$$
Now if $Y'$ is a finite set and $V\in\CC(Y',Y)$, the commutative diagram~(\ref{diagram}), for $X'=X$ and $U=\Delta_X$, becomes
\sbox{\matrice}{\begin{footnotesize}$\left(\!\!\begin{array}{cc}\Delta_X&\!\!\!\!\emptyset\\ \emptyset&\!\!\!\! V\end{array} \!\!\!\right)$\end{footnotesize}}
$$\xymatrix{
M'(X)\otimes_kM(Y)\ar[d]_{\Delta_X\otimes V}\ar[r]^-{\widehat{\psi}_{X,Y}}
&M''(X\sqcup Y)\quad \ar[d]^{\usebox{\matrice}} \ar[r]^-{=} & \; M''_Y(X) \ar[d]^{M''_V} \\
M'(X)\otimes_kM(Y')\ar[r]^-{\widehat{\psi}_{X,Y'}} &M''(X\sqcup Y') \quad \ar[r]^-{=} & \; M''_{Y'}(X)\mvirg\\
}$$
and it follows that the maps $\sur{\psi}_Y$ define a morphism of correspondence functors $\sur{\psi}:M\to \CH(M',M'')$.\par

Conversely, a morphism of correspondence functors $\xi:M\to \CH(M',M'')$ is determined by maps
$$\xi_Y:M(Y)\longrightarrow \CH(M',M'')(Y)=\Hom_{\CF_k}(M',M''_Y) \mvirg$$
for all finite sets $Y$.
Furthermore, for $m\in M(Y)$, the morphism $\xi_Y(m)$ is in turn determined by maps
$$\xi_Y(m)_X:M'(X)\longrightarrow M''_Y(X)=M''(X\sqcup Y)$$
for all finite sets $X$.
We claim that the family of maps
$$\xymatrix@R=1ex{
\mathring{\xi}_{X,Y}:M'(X)\otimes_k M(Y)\ar[r]&\;M''(X\sqcup Y)\\
{\phantom{\mathring{\xi}_{X,Y}:}} \hspace{5ex}m'\otimes m\hspace{5ex}\ar@{|->}[r]&\;\xi_Y(m)_X(m')
}
$$
defines a bilinear pairing $\mathring{\xi}:M'\times M\to M''$.
We must show that, for any finite sets $X,Y,X',Y'$ and any correspondences $U\in\CC(X',X)$ and $V\in\CC(Y',Y)$, the diagram
\sbox{\matrice}{\begin{footnotesize}$\left(\!\!\begin{array}{cc} U&\!\!\!\!\emptyset\\\emptyset&\!\!\!\!V\end{array}\!\!\!\right)$\end{footnotesize}}
$$\xymatrix{
M'(X)\otimes_kM(Y)\ar[d]_{U\otimes V}\ar[r]^-{\mathring{\xi}_{X,Y}}&M''(X\sqcup Y)\ar[d]^{\usebox{\matrice}}\\
M'(X')\otimes_kM(Y')\ar[r]^-{\mathring{\xi}_{X',Y'}}&M''(X'\sqcup Y')\\
}
$$
is commutative.
First observe that we have
\begin{eqnarray*}
\left(\begin{array}{cc} U&\emptyset\\\emptyset&V\end{array}\right) \mathring{\xi}_{X,Y}(m'\otimes m) &=&
\left(\begin{array}{cc} U&\emptyset\\\emptyset&V\end{array}\right) \xi_Y(m)_X(m') \\
&=&\left(\begin{array}{cc} \Delta_{X'}&\emptyset\\\emptyset&V\end{array}\right)
\left(\begin{array}{cc} U&\emptyset\\\emptyset&\Delta_Y\end{array}\right) \xi_Y(m)_X(m') \\
&=&\left(\begin{array}{cc} \Delta_{X'}&\emptyset\\\emptyset&V\end{array}\right) \xi_Y(m)_{X'}(Um') \mvirg
\end{eqnarray*}
because $\xi_Y(m)$ is a morphism of correspondence functors $M'\to M''_Y$.
Now the action of the correspondence $\left(\begin{array}{cc} \Delta_{X'}&\emptyset\\\emptyset&V\end{array}\right)$ is the composition with
$M''_V:M''_Y(X')\to M''_{Y'}(X')$,
which is in turn the action of~$V$ within the correspondence functor $\CH(M',M'')$.
Therefore we obtain
$$\left(\begin{array}{cc} \Delta_{X'}&\emptyset\\\emptyset&V\end{array}\right) \xi_Y(m)_X(Um')
= \xi_Y(Vm)_{X'}(Um') = \mathring{\xi}_{X,Y}(Um'\otimes Vm) \mvirg$$
using the fact that $\xi :M\to \CH(M',M'')$ is a morphism of correspondence functors.
This proves the claim.

By Theorem~\ref{pairings}, the pairing $\mathring{\xi}$ defines a morphism of correspondence functors
$$\check{\xi}:=\widetilde{\mathring{\xi}}: M'\otimes M \longrightarrow M'' \mpoint$$
Now it is straightforward to check that the maps $\psi\mapsto\sur{\psi}$ and $\xi\mapsto \check{\xi}$ are inverse isomorphisms
between $\Hom_{\CF_k}(M'\otimes M,M'')$ and $\Hom_{\CF_k}\!\big(M,\CH(M',M'')\big)$.
\endpf

In the case of a representable functor $k\CC(-,E)$, there is the following useful isomorphism.

\result{Proposition} \label{iso-repr} Let $E$ be a finite set and $N$ a correspondence functor over~$k$.
There is an isomorphism of correspondence functors $\CH\big(k\CC(-,E),N\big)\cong N_E$.
\fresult

\pf
Let $X$ be a finite set. By Yoneda's lemma, we get
$$\CH\big(k\CC(-,E),N\big)(X)=\Hom_{\CF_k}\big(k\CC(-,E),N_X\big)\cong N_X(E) \mpoint$$
Moreover, $N_X(E)=N(E\sqcup X)\cong N(X\sqcup E)=N_E(X)$.
It is straightforward to check that the resulting isomorphism 
$$\CH\big(k\CC(-,E),N\big)(X)\cong N_E(X)$$
is compatible with correspondences,
so that it yields an isomorphism of correspondence functors $\CH\big(k\CC(-,E),N\big)\cong N_E$.
\endpf

\result{Corollary} \label{Hom-constant-1} Let $N$ be a correspondence functor over~$k$.
There is an isomorphism of correspondence functors $\CH(\sou{k},N)\cong N$.
\fresult

\pf
Take $E=\emptyset$ in Proposition~\ref{iso-repr}. Then
$k\CC(-,\emptyset)\cong \sou{k}$ because $\CC(X,\emptyset)$ is the set of subsets of~$\emptyset$, which is a singleton.
On the other hand we clearly have $N_\emptyset=N$.
\endpf

Taking the constant functor~$\sou{k}$ in the second variable, we obtain a quite different result.

\result{Proposition} \label{Hom-constant-2} Let $M$ be a correspondence functor over~$k$.
There is an isomorphism of correspondence functors
$$\CH(M,\sou{k})\cong \sou{k}\otimes_k M(\emptyset)\dual$$
where $M(\emptyset)\dual=\Hom_k\big(M(\emptyset),k\big)$.
In particular, if $k$ is a field and if $M(\emptyset)$ is finite-dimensional, $\CH(M,\sou{k})$ is isomorphic to a direct sum of $\dim\!\big(M(\emptyset)\big)$ copies of~$\sou{k}$.
\fresult

\pf
Let $E$ be a finite set.
It is straightforward to see that $\sou{k}_E \cong \sou{k}$. Therefore
$$\CH(M,\sou{k})(E)\cong \Hom_{\CF_k}(M,\sou{k}_E)\cong \Hom_{\CF_k}(M,\sou{k})\cong \Hom_k\big(M(\emptyset),k\big)=M(\emptyset)\dual \mpoint$$
It is then easy to check that the action of a correspondence $U\in \CC(F,E)$ yields the identity endomorphism of $M(\emptyset)\dual$,
so that we get the constant functor tensored with~$M(\emptyset)\dual$.
If $k$ is a field and if $M(\emptyset)$ is finite-dimensional,
it follows that $\CH(M,\sou{k})$ is isomorphic to a direct sum of copies of~$\sou{k}$, their number being
$\dim\!\big(M(\emptyset)\dual\big)=\dim\!\big(M(\emptyset)\big)$.
\endpf

In the same vein as in Theorem~\ref{finite-gen}, we now consider bounded type and finite generation.

\result{Theorem} \label{finite-gen-hom} Let $M$ and $N$ be correspondence functors over~$k$.
\begin{enumerate}
\item Let $E$ and $F$ be finite sets. If $M$ is generated by~$M(E)$, then $M_F$ is generated by~$M_F(E)$.
Therefore, if $M$ has bounded type, so has~$M_F$.
If $M$ is finitely generated, so is~$M_F$.
\item Assume that the ring~$k$ is noetherian.
If $M$ is finitely generated and if $N$ has bounded type (respectively is finitely generated),
then $\CH(M,N)$ has bounded type (respectively is finitely generated).
\end{enumerate}
\fresult

\pf (a) By assumption, $M(X)=k\CC(X,E)M(E)$ for each finite set $X$.
Replacing $X$ by $X\sqcup F$ gives
$$M(X\sqcup F)=k\CC(X\sqcup F,E)M(E)\mpoint$$ 
Therefore $M(X\sqcup F)$ is $k$-linearly generated by the elements $\displaystyle \binom{V}{W}(m)$, where $V\in\CC(X,E)$, $W\in\CC(F,E)$, and $m\in M(E)$,
because any correspondence in $\CC(X\sqcup F,E)$ can be written $\displaystyle \binom{V}{W}$.
But we have
$$\binom{V}{W}(m)=\left(\begin{array}{cc}V&\emptyset\\\emptyset&\Delta_F\end{array}\right)\binom{\Delta_E}{W}(m)$$
and this is the image of ${\displaystyle \binom{\Delta_E}{W}}(m)\in M_F(E)$ by the correspondence $V$ within the functor~$M_F$.
It follows that
$$\binom{V}{W}(m)\in k\CC(X,E)M_F(E) \mpoint$$
Therefore $M(X\sqcup F)$ is $k$-linearly generated by $k\CC(X,E)M_F(E)$, that is,
$$M_F(X)=k\CC(X,E)M_F(E) \mvirg$$
as was to be shown.
The other two assertions follow immediately.\mpn
 
(b) Since $M$ is finitely generated, there is a finite subset $A\subseteq M(E)$ such that $M=k\CC(-,E)A$,
so $M$ is isomorphic a quotient of the finite direct sum $\dirsum{a\in A} k\CC(-,E)$ of representable functors.
Since $\CH(-,N)$ is exact by Lemma~\ref{internal-hom}, we deduce an embedding
$$\xymatrix{
\CH(M,N) \; \ar@{^{(}->}[r]&\dirsum{a\in A} \CH\big(k\CC(-,E),N\big)\cong\dirsum{a\in A} N_{E}\mvirg
}
$$
using also Proposition~\ref{iso-repr}.
If $N$ has bounded type, then $N_E$ has bounded type, by~(a), so the finite direct sum $\dirsum{a\in A} N_E$ also has bounded type.
Therefore $\CH(M,N)$ is isomorphic to a subfunctor of a functor of bounded type.
Since $k$ is noetherian, this implies that $\CH(M,N)$ has bounded type, by Corollary~11.5 in~\cite{BT2}.
The same argument with ``bounded type'' replaced by ``finitely generated'' goes through, completing the proof.
\endpf


\bigbreak

\section{Algebra functors}\label{Section:Algebra}

\medskip
\noindent
For any finite lattice~$T$, the functor $F_T$ has more structure, namely it is a commutative algebra in the tensor category~$\CF_k$ of all correspondence functors. This section is devoted to a closer analysis of this additional structure.
By a $k$-algebra, we always mean an associative $k$-algebra with an identity element.

\result{Definition} An {\em algebra correspondence functor over~$k$} is a correspondence functor~$A$
with values in the category $k\text{-\!}\Alg$ of $k$-algebras satisfying the following two conditions.
For any finite sets $X$ and $Y$, and for any correspondence $U\in\CC(Y,X)$, the diagrams
$$\xymatrix{
A(X)\otimes_k A(X)\ar[r]^-{\mu_X}\ar[d]_-{U\otimes U} &A(X)\ar[d]^-{U}\\
A(Y)\otimes_k A(Y)\ar[r]^-{\mu_Y}&A(Y)\\
}
\qquad\text{and}\qquad
\xymatrix{
k \ar[r]^-{\varepsilon_X}\ar[d]_-{\Id} &A(X)\ar[d]^-{U}\\
k \ar[r]^-{\varepsilon_Y} &A(Y)
}
$$
are commutative, where $\mu_X:A(X)\otimes_k A(X) \to A(X)$ denotes the multiplication map of the $k$-algebra $A(X)$
and $\varepsilon_X: k\to A(X)$ denotes the map $\lambda\mapsto \lambda\cdot 1_{A(X)}$.

An algebra correspondence functor $A$ is called {\em commutative} if $A(X)$ is a commutative $k$-algebra for any finite set~$X$.
\fresult

The commutativity of all the diagrams in the definition can be interpreted in two different ways~:
\begin{enumerate}
\item The action of any correspondence $U\in\CC(Y,X)$ is a map of $k$-algebras $A(X) \to A(Y)$,
mapping $1_{A(X)}$ to~$1_{A(Y)}$.
\item The family of multiplication maps $\mu_X$ defines a morphism of correspondence functor $\mu:A\otimes A \to A$
and the family of maps $\varepsilon_X$ defines a morphism of correspondence functors
$\varepsilon:\sou{k} \longrightarrow A$,
where $\sou{k}$ denotes the constant functor of Example~\ref{constant}.
In other words, the triple $(A,\mu,\varepsilon)$ is an algebra in the tensor category $\CF_k$.
\end{enumerate}

\result{Lemma} \label{mu}
Let $A$ be an algebra correspondence functor, let $\mu:A\otimes A\to A$ be the multiplication map,
let $\widehat\mu :A\times A \to A$ be the associated pairing (see Theorem~\ref{pairings}),
and let $X$ and $Y$ be finite sets.
\begin{enumerate}
\item The map
$$\widehat\mu_{X,Y} : A(X)\otimes_k A(Y) \longrightarrow A(X\sqcup Y) \mpoint$$
is obtained as the composite of the action of $\binom{\Delta_X}{\emptyset}\otimes \binom{\emptyset}{\Delta_Y}$
and the multiplication $\mu_{X\sqcup Y}$.
\item If $A$ is commutative, then $\widehat\mu_{X,Y}$ is an algebra homomorphism.
\item The composite $(\Delta_X,\Delta_X)\,\widehat\mu_{X,X}$ is equal to~$\mu_X$.
\end{enumerate}
\fresult

\pf
(a) This follows from the proof of Theorem~\ref{pairings}.\mpn

(b) The action of $\binom{\Delta_X}{\emptyset}\otimes \binom{\emptyset}{\Delta_Y}$ is an algebra homomorphism and,
whenever $A$ is commutative, so is the multiplication $\mu_{X\sqcup Y}$.\mpn

(c) Using (a), we obtain
\begin{eqnarray*}
(\Delta_X,\Delta_X) \, \widehat\mu_{X,X}
&=& (\Delta_X,\Delta_X)\, \mu_{X\sqcup Y} \Big(\binom{\Delta_X}{\emptyset}\otimes \binom{\emptyset}{\Delta_Y} \Big) \\
&=& \mu_X\Big((\Delta_X,\Delta_X)\otimes (\Delta_X,\Delta_X)\Big)
\Big(\binom{\Delta_X}{\emptyset}\otimes \binom{\emptyset}{\Delta_Y} \Big) \\
&=& \mu_X (\Delta_X\otimes\Delta_X) =\mu_X \mpoint
\end{eqnarray*}
\endpf

We consider now the functor $F_T$ of Example~\ref{F_T}, associated to a finite lattice~$T$.

\result{Proposition} \label{characterize}
Let $T$ be a finite lattice and let $\hat 0$ be its least element.
\begin{enumerate}
\item {\bf Commutative algebra.} $F_T$ is a commutative algebra correspondence functor, with respect to the multiplication maps
$$\mu_X: F_T(X)\otimes_k F_T(X) \to F_T(X) \,, \qquad \varphi\otimes\psi \mapsto \varphi\vee\psi \mvirg$$
where $X$ is a finite set and $\varphi,\psi: X\to T$ are maps.
Here $\varphi\vee\psi:X\to T$ denotes the map defined by
$(\varphi\vee\psi)(x)=\varphi(x)\vee\psi(x)$ for every $x\in X$.
The identity element of $F_T(X)$ is the constant map onto~$\hat 0$.
\item {\bf Exponential property.} $F_T(\emptyset)\cong k$ and, for any finite sets $X$ and~$Y$, the bilinear pairing
$$\widehat{\mu}_{X,Y}:F_T(X) \otimes_k F_T(Y) \to F_T(X\sqcup Y)$$
associated with~$\mu$ is an isomorphism of $k$-modules.\par

More precisely, for any maps $\varphi:X\to T$ and $\psi:Y\to T$, the element $\widehat{\mu}_{X,Y}(\varphi\otimes\psi)\in F_T(X\sqcup Y)$ is the function $X\sqcup Y \to T$ equal to $\varphi$ on~$X$ and to $\psi$ on~$Y$,
providing a bijection between the canonical bases of~$F_T(X) \otimes_k F_T(Y)$ and~$F_T(X\sqcup Y)$.
\item {\bf Splitting property.} If $\bullet$ denotes a set of cardinality one, then $F_T(\bullet)$ is isomorphic to a finite direct product $k\times k\times \ldots \times k$ as $k$-algebras.
\end{enumerate}
\fresult

\pf (a) We claim that the map
$$\upsilon :T\times T\to T \,,\quad\upsilon(a,b)=a\vee b \mvirg$$
is a join-morphism.
Given a subset $C\subseteq T\times T$, let $A\subseteq T$ be its first projection and $B\subseteq T$ its second projection.
Then
\begin{eqnarray*}
\upsilon\big(\mbigvee{(a,b)\in C}(a,b)\big)&=&\upsilon\Big(\mbigvee{(a,b)\in C}\big((a,\hat0)\vee(\hat0,b)\big)\Big)\\
&=&\upsilon\Big(\big(\mbigvee{a\in A}(a,\hat0)\big) \vee \big(\mbigvee{b\in B}(\hat0,b) \big) \Big)\\
&=&\upsilon(\big(\mbigvee{a\in A}a\,,\mbigvee{b\in B}b\big)\\
&=&(\mbigvee{a\in A} a) \vee (\mbigvee{b\in B}b) \\
&=&\mbigvee{(a,b)\in C}(a\vee b) \\
&=&\mbigvee{(a,b)\in C}\upsilon(a,b) \mvirg
\end{eqnarray*}
proving the claim.
By Example~\ref{F_T}, the map $\upsilon$ induces a morphism $F_{T\times T}\to F_T$.
Since $F_{T\times T}\cong F_T\otimes F_T$ by Theorem~\ref{tensor lattice},
we obtain a morphism $\mu:F_T\otimes F_T\to F_T$. For any finite set~$X$,
the map
$$\mu_X:F_T(X)\otimes_k F_T(X) \to F_T(X)$$
is easily seen to be the map of the statement.
Clearly $\mu_X$ is associative and commutative and the constant map onto~$\hat0$ is an identity element.
We obtain in this way an algebra correspondence functor $F_T$ because $\mu:F_T\otimes F_T\to F_T$ is a morphism of functors,
and so is $\varepsilon:\sou{k} \to F_T$.\mpn

(b) One checks easily that $\binom{\Delta_X}{\emptyset}\varphi$ is the function from $X\sqcup Y$ to $T$ equal to $\varphi$ on $X$, and to $\hat0$ on $Y$. Similarly $\binom{\emptyset}{\Delta_Y}\psi$ is the map equal to $\hat0$ on $X$ and to $\psi$ on $Y$. Thus $\Big(\binom{\Delta_X}{\emptyset}\varphi\Big)\vee\Big(\binom{\emptyset}{\Delta_Y}\psi\Big)$ is the map equal to $\varphi$ on $X$ and to $\psi$ on $Y$.\mpn

(c) Since $\bullet$ has cardinality one, $F_T(\bullet)$ is a free $k$-module with basis
$$\{g_t\mid t\in T\}$$
where $g_t:\bullet \to T$ is defined by $g_t(\bullet)=t$ (with $\bullet$ being also the unique element of the set~$\bullet$).
Moreover, on restriction to this basis, the multiplication map corresponds to the map~$\upsilon$ of the beginning of the proof,
namely $g_t g_{t'}=g_{t\vee t'}$.
In that case, there is a standard procedure for finding another $k$-basis consisting of orthogonal primitive idempotents~$f_t \in T$
whose sum is the identity element (namely~$\hat0$). These primitive idempotents are defined by
$$f_t=\sum_{\substack{s\in T \\ s\geq t}} \chi(t,s) g_s \mvirg$$
where $\chi(t,s)$ denotes the M\"obius function of the poset~$T$ (see the appendix in~\cite{BT1}).
Then we obtain isomorphisms of $k$-algebras
$$F_T(\bullet) \;\cong\; \prod_{t\in T} \,k{\cdot} f_t \;\cong\; k\times k\times \ldots \times k \mvirg$$
as required.
\endpf

Our next main result asserts that, with a small assumption on~$k$, the converse of Proposition~\ref{characterize} holds.

\result{Theorem} \label{converse} Assume that $k$ does not contain a nontrivial idempotent (i.e. ${\rm Spec}(k)$ is connected).
Let $A$ be a correspondence functor over~$k$ and suppose that $A$ has the following three properties~:
\begin{enumerate}
\item {\bf Commutative algebra.} $A$ is a commutative algebra correspondence functor (with multiplication written~$\mu$).
\item {\bf Exponential property.} For any finite sets $X$ and~$Y$, the associated bilinear pairing
$$\widehat{\mu}_{X,Y}:A(X) \otimes_k A(Y) \to A(X\sqcup Y)$$
is an isomorphism of $k$-modules (hence an isomorphism of $k$-algebras, by commutativity of~$A$).
Moreover, $A(\emptyset)\cong k$.
\item {\bf Splitting property.} If $\bullet$ denotes a set of cardinality one,
then $A(\bullet)$ is isomorphic to a finite direct product $k\times k\times \ldots \times k$ as $k$-algebras.

\end{enumerate}
Then there exists a finite lattice $T$ such that $A\cong F_T$ (isomorphism of algebra correspondence functors).
\fresult

\bigskip
We need a preliminary lemma.

\result{Lemma} \label{comultiplication}
With the assumptions above, define a comultiplication
$$\delta_\bullet : A(\bullet) \longrightarrow A(\bullet) \otimes_k A(\bullet)$$
as the composition of the action of $\binom{\Delta_\bullet}{\Delta_\bullet}:A(\bullet) \to A(\bullet\sqcup\bullet)$ and the isomorphism
$$\widehat{\mu}_{\bullet,\bullet}^{-1}: A(\bullet\sqcup\bullet) \longrightarrow A(\bullet) \otimes_k A(\bullet) \mpoint$$
\begin{enumerate}
\item $\delta_\bullet$ is an algebra homomorphism.
\item $\delta_\bullet$ is coassociative and cocommutative.
\item $\mu_\bullet \delta_\bullet = \Id_{A(\bullet)}$.
\item The map $\eta_\bullet: A(\bullet)\to A(\emptyset)=k$ induced by the action of the empty correspondence $\emptyset\in\CC(\emptyset,\bullet)$
is a counit for the comultiplication~$\delta_\bullet$ and is an algebra homomorphism.
\end{enumerate}
\fresult

\pf
(a) The action of $\binom{\Delta_\bullet}{\Delta_\bullet}:A(\bullet) \to A(\bullet\sqcup\bullet)$ is an algebra homomorphism. So is the map ${\widehat\mu}_{\bullet,\bullet}$ by Lemma~\ref{mu},
hence so is the composite ${\widehat\mu}_{\bullet,\bullet}^{-1} \binom{\Delta_\bullet}{\Delta_\bullet}$.\mpn

(b) This is left as an exercise for the reader. For the cocommutativity, use both the twist of $A(\bullet)\otimes_k A(\bullet)$ 
and the action of the correspondence
$\left(\begin{array}{cc}\emptyset &\Delta_\bullet \\ \Delta_\bullet & \emptyset \end{array}\right)$.\mpn

(c) We compute
$$\mu_\bullet \delta_\bullet
=(\Delta_\bullet,\Delta_\bullet){\widehat\mu}_{\bullet,\bullet} {\widehat\mu}_{\bullet,\bullet}^{-1}
\binom{\Delta_\bullet}{\Delta_\bullet}=(\Delta_\bullet,\Delta_\bullet)\binom{\Delta_\bullet}{\Delta_\bullet}
=\Id_{A(\bullet)} \mpoint$$

(d) A special case of diagram~(\ref{diagram}) yields the commutative diagram
$$\xymatrix{
& A(\bullet)\otimes_k A(\bullet)\ar[d]_{\eta_\bullet\otimes \Delta_\bullet} \ar[r]^-{{\widehat\mu}_{\bullet,\bullet}}
&A(\bullet\sqcup \bullet)\ar[d]^{(\emptyset,\Delta_\bullet)}\\
A(\bullet) \ar[r]^-{\sim} &A(\emptyset)\otimes_k A(\bullet)\ar[r]^-{{\widehat\mu}_{\emptyset,\bullet}}
&A(\emptyset\sqcup \bullet) \ar[r]^-{=} &A(\bullet)
}
$$
and it easy to check that the bottom composite is the identity~$\Id_{A(\bullet)}$.
We deduce the commutative diagram
$$\xymatrix{
A(\bullet) \ar[r]^-{\binom{\Delta_\bullet}{\Delta_\bullet}} )\ar[d]_{=}
& A(\bullet\sqcup \bullet)\ar[d]_{(\emptyset,\Delta_\bullet)} \ar[r]^-{{\widehat\mu}_{\bullet,\bullet}^{-1}}
&A(\bullet)\otimes_k A(\bullet) \ar[d]^{\eta_\bullet\otimes \Delta_\bullet}\\
A(\bullet) \ar[r]^-{=} & A(\emptyset\sqcup \bullet) \ar[r]^-{{\widehat\mu}_{\emptyset,\bullet}^{-1}}
&A(\emptyset)\otimes_k A(\bullet) \ar[r]^-{\sim} &A(\bullet)
}
$$
and since the composite of the first row is equal to~$\delta_\bullet$,
we obtain that $(\eta_\bullet\otimes \Delta_\bullet)\delta_\bullet$ is the identity~$\Id_{A(\bullet)}$, as required for a counit.
\endpf

\bigskip\noindent{\bf Proof of Theorem~\ref{converse}~: }
By the splitting property, the $k$-algebra $A(\bullet)$ contains a $k$-basis
$$\{ f_t\mid t\in T\}$$
consisting of orthogonal idempotents such that $\sum_{t\in T} f_t=1_{A(\bullet)}$, where $T$ is a finite index set.
Our first aim is to show that $T$ has a lattice structure.\par

Now $\{f_a\otimes f_b \mid a,b\in T\}$ is a $k$-basis of orthogonal idempotents of~$A(\bullet) \otimes_k A(\bullet)$,
so we can write
$$\delta_\bullet(f_t) = \sum_{a,b\in T} \lambda_{a,b}\, f_a\otimes f_b \,,\qquad \lambda_{a,b}\in k \mpoint$$
Since $\delta_\bullet : A(\bullet) \longrightarrow A(\bullet) \otimes_k A(\bullet)$ is an algebra homomorphism
by Lemma~\ref{comultiplication}, $\delta_\bullet(f_t)$ is an idempotent of~$A(\bullet) \otimes_k A(\bullet)$.
Therefore $\lambda_{a,b}$ is an idempotent of~$k$, hence $\lambda_{a,b}\in\{0,1\}$ by our assumption on~$k$, and
$$\delta_\bullet(f_t) = \sum_{(a,b)\in B_t} f_a\otimes f_b$$
where $B_t$ is a subset of $T\times T$.
Since the idempotents $f_t$ are orthogonal and sum to~$1_{A(\bullet)}$, the idempotents $\delta_\bullet(f_t)$ are orthogonal
and sum to $1_{A(\bullet)}\otimes 1_{A(\bullet)}$.
This shows that $B_s \cap B_t = \emptyset$ for any pair $(s,t)$ of distinct elements of~$T$, and that $\displaystyle\bigsqcup_{t\in T} B_t = T \times T$.
Hence we can define a operation $\meet$ on~$T$ by
$$a\meet b=t \;\Longleftrightarrow\; (a,b) \in B_t \mpoint$$
In other words, for any $t\in T$,
$$\delta_\bullet(f_t)=\sum_{\substack{a,b\in T \\ a\meet b=t}} f_a\otimes f_b \mpoint$$
Since $\delta_\bullet$ is coassociative by Lemma~\ref{comultiplication}, this operation~$\meet$ is associative.
Similarly, since $\delta_\bullet$ is cocommutative, $\meet$ is commutative.
Finally, by Lemma~\ref{comultiplication} again, for any $t\in T$, we have
$$f_t=\mu_\bullet \delta_\bullet(f_t)= \sum_{\substack{a,b\in T \\ a\meet b=t}} f_af_b
=\sum_{\substack{a\in T \\ a\meet a=t}} f_a \mvirg$$
because $f_af_b=0$ if $a\neq b$.
It follows that $t\meet t=t$ for any $t\in T$.
Hence $T$ is a commutative idempotent semigroup.
Equivalently, if we define a relation~$\leq$ on~$T$ by
$$a \leq b \;\Longleftrightarrow\; a\meet b = a \mvirg$$
we get an order relation on~$T$, and any pair $\{a,b\}$ of elements of~$T$ has a greatest lower bound $a\meet b$.
Thus $T$ is a meet semilattice.\par

Since the counit $\eta_\bullet: A(\bullet)\to k$ is an algebra homomorphism by Lemma~\ref{comultiplication}
and since the only idempotents of~$k$ are 0 and~1,
we have $\eta_\bullet(f_t)\in \{0,1\}$ for any $t\in T$.
Moreover
$$1=\eta_\bullet(1_{A(\bullet)}) = \eta_\bullet\big(\sum_{t\in T} f_t\big) = \sum_{t\in T} \eta_\bullet(f_t) \mvirg$$
so there exists $u \in T$ such that $\eta_\bullet(f_u) = 1$.
This element $u$ is unique because if $u, v\in T$ are such that $\eta_\bullet(f_u) =\eta_\bullet(f_v) = 1$,
then $\eta_\bullet(f_uf_v) = 1$, so $f_uf_v \neq0$, hence $u = v$.
Moreover, since $\eta_\bullet$ is a counit, we have, for any $t\in T$,
$$f_t=(\eta_\bullet\otimes\Id)\delta_\bullet(f_t)=\sum_{\substack{a,b\in T \\ a\meet b = t}} \eta_\bullet(f_a)f_b
=\sum_{\substack{b\in T \\ u\meet b = t}} f_b$$
and it follows that $u\meet t=t$, hence $t\leq u$.
Therefore $T$ is a meet semilattice with a greatest element~$u$, so it is a lattice.
We write $\vee$ for its join operation and $\hat0$ for its least element.\par

Now, for any $t\in T$, we define
$$g_t=\sum_{\substack{s\in T \\ s\geq t}} f_s \mpoint$$
The elements $\{g_t\mid t\in T\}$ form another basis of~$A(\bullet)$, because the transition matrix is unitriangular.
Now for any $t,t'\in T$, we have
$$g_t g_{t'}=\sum_{\substack{s\geq t \\ s'\geq t'}} f_s f_{s'} = \sum_{s\geq t\vee t'} f_s = g_{t\vee t'}$$
and also
$$g_{\hat0}=\sum_{\substack{s\in T \\ s\geq \hat0}} f_s =\sum_{s\in T}f_s=1_{A(\bullet)} \mpoint$$
This means that the map $t\mapsto g_t$, from $T$ to~$A(\bullet)$, induces an algebra homomorphism
$$\theta_\bullet: F_T(\bullet)=kT^\bullet \longrightarrow A(\bullet) \,,\qquad \varphi \mapsto g_{\varphi(\bullet)} \,,\; \forall\, \varphi\in T^\bullet \mvirg$$
and this is an isomorphism because it maps a basis to a basis.\par

Whenever we have a disjoint union $X=W\sqcup Z$, there is a bilinear pairing
$$\widehat{\mu}_{W,Z}:A(W) \otimes_k A(Z) \to A(X) \mvirg$$
which is an isomorphism by the exponential property.
Decomposing $X$ as a union of singletons, we obtain by induction an isomorphism of $k$-algebras
$$\widehat\mu: \; \bigotimes_{x\in X} A(\bullet) \longrightarrow A(X) \mpoint$$
By part~(a) of Lemma~\ref{mu}, the isomorphism $\widehat\mu$ maps $\displaystyle\bigotimes_x a_x$ to the product $\displaystyle\prod_{x\in X} C_x a_x$, where $C_x:=\{ (x,\bullet) \} \subseteq X\times \bullet$\,.
We then obtain a sequence of isomorphisms of $k$-algebras
$$\xymatrix{
F_T(X) \ar[r]_-{\widehat\mu^{-1}}^-{\sim} & \displaystyle \;\bigotimes_{x\in X} F_T(\bullet)\; \ar[r]_-{\displaystyle\otimes_x\, \theta_{\bullet}}^-{\sim}
& \displaystyle\;\bigotimes_{x\in X} A(\bullet)\; \ar[r]_-{\widehat\mu}^-{\sim} & A(X) \mpoint
}
$$
It is easy to check that the first isomorphism $\mu^{-1}$ maps a function $\varphi:X\to T$ to the element
$\displaystyle\bigotimes_{x\in X} h_{\varphi(x)}$, where $h_{\varphi(x)}: \bullet \to T$ is defined by $h_{\varphi(x)}(\bullet)=\varphi(x)$.
The second isomorphism maps $\displaystyle\bigotimes_{x\in X} h_{\varphi(x)}$ to $\displaystyle\bigotimes_{x\in X} g_{\varphi(x)}$,
which in turn maps to $\displaystyle\prod_{x\in X} C_x g_{\varphi(x)}$ via the third isomorphism~$\widehat\mu$.
Therefore we get the composite isomorphism of $k$-algebras
$$\lambda_X:F_T(X) \longrightarrow A(X) \,,\qquad \lambda_X(\varphi)=\prod_{x\in X} C_x g_{\varphi(x)} \;, \;\forall \varphi\in T^X \mpoint$$
We are going to show that the isomorphisms $\lambda_X$ are compatible with the action of correspondences, but we first need a lemma.

\result{Lemma} \label{product=union}
Let $Y$ be a finite set and let $W$ and~$Z$ be subsets of $Y\times\bullet$.
Let $t \in T$ and let $g_t\in A(\bullet)$ as defined above. Then, in the $k$-algebra $A(Y)$, we have an equality
$$(W g_t)(Z g_t) = (W\cup Z) g_t \mpoint$$
\fresult

\pf
By (\ref{diagram}), there is a commutative diagram
\sbox{\matrice}{\begin{footnotesize}$\left(\!\!\begin{array}{cc}W &\!\!\!\!\emptyset\\ \emptyset&\!\!\!\! Z\end{array}\!\!\right)$\end{footnotesize}}
$$\xymatrix{A(\bullet) \ar[r]^-{\delta_\bullet} & A(\bullet)\otimes_k A(\bullet) \ar[r]^-{{\widehat\mu}_{\bullet,\bullet}} \ar[d]_-{W\otimes Z}
& A(\bullet\sqcup\bullet) \;\;\; \ar[d]^-{\usebox{\matrice}} \\
& A(Y)\otimes_k A(Y) \ar[r]_-{{\widehat\mu}_{Y,Y}} & A(Y\sqcup Y) \;\;\; \ar[r]_-{(\Delta_Y,\Delta_Y)} & \;\;A(Y)
}
$$
and we compute the image of~$g_t$ using both paths.
For the top path, the definition of~$\delta_\bullet$ gives
$${\widehat\mu}_{\bullet,\bullet} \,\delta_\bullet(g_t)=\binom{\Delta_\bullet}{\Delta_\bullet} g_t$$
and since $(\Delta_Y,\Delta_Y) \usebox{\matrice} \binom{\Delta_\bullet}{\Delta_\bullet} = W\cup Z$, we obtain the element
$$(W\cup Z)g_t \in A(Y) \mpoint$$
For the bottom path, we first have
$$\delta_\bullet(g_t)=\sum_{s\geq t} \delta_\bullet(f_s) =\sum_{s\geq t} \sum_{\substack {a,b\in T \\ a\meet b=s}} f_a\otimes f_b
=\sum_{\substack {a,b\in T \\ a\meet b\geq t}} f_a\otimes f_b = \big(\sum_{a\geq t} f_a \big) \otimes \big(\sum_{b\geq t} f_b\big)
= g_t\otimes g_t $$
and this is mapped to $Wg_t\otimes Zg_t$ by the vertical map $W\otimes Z$.
Since the composition of the two bottom maps is $(\Delta_Y,\Delta_Y) {\widehat\mu}_{Y,Y}= \mu_Y$, we obtain
$$\mu_Y(Wg_t\otimes Zg_t)=(Wg_t) (Zg_t) \in A(Y) \mpoint$$
This completes the proof of the lemma.
\endpf

Now we return to the proof of the theorem and we let $U\subseteq Y\times X$ and $\varphi\in T^X$.
Since $(U\varphi)(y)=\displaystyle\bigvee_{\substack {x\in X \\ (y,x)\in U}} \varphi(x)$ and since $g_{a\vee b}=g_ag_b\,, \;\forall a,b\in T$, we obtain
\begin{eqnarray*}
\lambda_Y (U\varphi) &=& \prod_{y\in Y} C_y \, g_{(U\varphi)(y)} = \prod_{y\in Y} C_y  \prod_{\substack {x\in X \\ (y,x)\in U}}  g_{\varphi(x)} \\
&=& \prod_{x\in X} \Big( \prod_{\substack {y\in Y \\ (y,x)\in U}} C_y \, g_{\varphi(x)} \Big) 
= \prod_{x\in X} \Big( \bigsqcup_{\substack {y\in Y \\ (y,x)\in U}} C_y\Big) g_{\varphi(x)}
\end{eqnarray*}
using Lemma~\ref{product=union}. Therefore
$$\lambda_Y (U\varphi) = \prod_{x\in X} \Big( \bigsqcup_{\substack {y\in Y \\ (y,x)\in U}} C_y\Big) g_{\varphi(x)}
=\prod_{x\in X} ( UC_x ) g_{\varphi(x)} = U \prod_{x\in X} C_x  g_{\varphi(x)} 
= U \lambda_X(\varphi) \mpoint$$
This proves that $\lambda:F_T\to A$ is an isomorphism of algebra correspondence functors.
\endpf

\result{Remark} \label{uniqueness}
{\rm
In Theorem~\ref{converse}, the algebra structure on $A(X)$ is uniquely determined from that of~$A(\bullet)$.
More precisely, if $A$ and $A'$ are correspondence functors satisfying the three properties of Theorem~\ref{converse} and if $A(\bullet)\cong A'(\bullet)$,
then $A\cong A'$. This follows from the proof, because the lattice $T$ is uniquely determined by the algebra structure of $A(\bullet)$.
}
\fresult

\result{Remark} \label{remark}
{\rm
The assumption on~$k$ is necessary in Theorem~\ref{converse}.
Let $k=k_1\times k_2$ be the direct product of two nontrivial rings
and, for $i=1,2$, let $F_{T_i}^{(k_i)}$ be the functor over~$k_i$ associated to a lattice~$T_i$,
but viewed as a functor over~$k$ (with the other factor $k_{3-i}$ acting by zero).
It can be shown that, if $T_1$ and $T_2$ have the same cardinality but are not isomorphic,
then $F_{T_1}^{(k_1)} \times F_{T_2}^{(k_2)}$ satisfies the three assumptions of Theorem~\ref{converse} but is not isomorphic to~$F_T$ for any lattice~$T$.
}
\fresult


\bigbreak

\section{Examples}\label{Section:Examples}

\medskip
\noindent
Finding the decomposition of tensor products is not straightforward, due in particular to fast increasing dimensions.
We only give here a few small examples, based on~\cite{BT3} and~\cite{BT4}.
We refer to those two papers for details.
For simplicity, we assume that $k$ is a field.\par

For any $n\in\N$, we let $[n]=\{1,2,\ldots,n\}$ and $\sou{n}=\{0\} \sqcup [n]$.
Then $\sou{n}$ is a totally ordered lattice and $[n]$ is its subset of irreducible elements.
There is a simple correspondence functor $\S_n$ introduced in Section~11 of~\cite{BT3}.
It appears as a direct summand of the functor~$F_{\sou{n}}$ associated to the lattice~$\sou{n}$.

\result{Example} \label{zero}
{\rm
When $n=0$, then $F_{\sou0}=\S_0$ is the constant functor~$\sou{k}$ of Example~\ref{constant}.
For any correspondence functor~$M$, we have $\S_0 \otimes M \cong M$, by Proposition~\ref{tensor-product}.
}
\fresult

\result{Example} \label{one}
{\rm
When $n=1$, we have $F_{\sou1}\cong \S_0 \oplus \S_1$, by Theorem~11.6 of~\cite{BT3}.
Moreover, by Theorem~\ref{tensor lattice}, we know that
$$F_{\sou1}\otimes F_{\sou1} \cong F_{\sou1\times\sou1} =F_{\lozenge} \mvirg$$
where $\lozenge:=\sou1\times\sou1$ denotes the lattice of subsets of a set of cardinality~2.
Applying Example~\ref{zero}, we obtain
\begin{eqnarray*}
F_\lozenge &\cong& F_{\sou1}\otimes F_{\sou1} \\
&\cong& (\S_0 \oplus \S_1)\otimes(\S_0 \oplus \S_1) \\
&\cong& \S_0 \oplus 2 \S_1 \oplus (\S_1\otimes \S_1)
\end{eqnarray*}
On the other hand, by Example~8.7 of~\cite{BT4}, there is a direct sum decomposition
\def\bcirc{{\scriptscriptstyle\circ\circ}}
$$F_\lozenge\cong \S_0\oplus 3\S_1\oplus 2\S_2\oplus \S_{\bcirc}\mvirg$$
where $\S_{\bcirc}$ is the fundamental functor associated to the poset $\bcirc$ of cardinality~2 ordered by the equality relation.\par

We now have two expressions for $F_\lozenge$ and we apply the Krull-Remak-Schmidt theorem, which holds by Proposition~6.6 in~\cite{BT2}.
It follows that
$$\S_1\otimes \S_1 \cong \S_1 \oplus 2 \S_2 \oplus \S_{\bcirc} \mpoint$$
}
\fresult

\result{Example} \label{two}
{\rm 
When $n=2$, we have $F_{\sou2}\cong \S_0 \oplus 2\S_1\oplus \S_2$, by Theorem~11.6 of~\cite{BT3}.
Moreover, by Theorem~\ref{tensor lattice}, we know that
$$F_{\sou1}\otimes F_{\sou2} \cong F_{\sou1\times\sou2} =F_P \mvirg$$
where $P:=\sou1\times\sou2$.
Applying the previous two examples, we obtain
\def\bcirc{{\scriptscriptstyle\circ\circ}}
\begin{eqnarray*}
F_P &\cong& F_{\sou1}\otimes F_{\sou2} \\
&\cong& (\S_0 \oplus \S_1)\otimes(\S_0 \oplus 2\S_1 \oplus \S_2) \\
&\cong& \S_0 \oplus 3 \S_1 \oplus \S_2 \oplus 2(\S_1\otimes \S_1) \oplus (\S_1\otimes \S_2) \\
&\cong& \S_0 \oplus 3 \S_1 \oplus \S_2 \oplus 2\S_1\oplus 4\S_2 \oplus 2\S_{\bcirc}  \oplus (\S_1\otimes \S_2) \\
&\cong& \S_0 \oplus 5 \S_1 \oplus 5\S_2 \oplus 2\S_{\bcirc}  \oplus (\S_1\otimes \S_2)
\end{eqnarray*}
On the other hand, by Example~8.11 of~\cite{BT4}, there is a direct sum decomposition
\def\posetVop{^{\mathop{\mathop{\scriptscriptstyle\circ\;\circ}^{\scriptscriptstyle/\backslash}\limits}^{\raisebox{-1ex}{$\scriptscriptstyle\circ$}}\limits}}
\def\posetV{^{\,\mathop{\mathop{\scriptscriptstyle\circ\;\circ}_{\scriptscriptstyle\backslash/}\limits}_{^\circ}\limits}}
$$F_P \cong \S_0\oplus 5\S_1\oplus 7\S_2\oplus3\S_3\oplus 3\S_{\bcirc}\oplus \S_{\posetV}\oplus \S_{\posetVop}\oplus U\mvirg$$
where $\S_{\posetV}$ (respectively $\S_{\posetVop}$) denotes the fundamental functor associated to the poset indicated as a subscript,
and where $U$ is an indecomposable projective functor of Loewy length~3 described in Example~8.11 of~\cite{BT4}.\par

We now have two expressions for $F_P$ and it follows from the Krull-Remak-Schmidt theorem that
$$\S_1\otimes \S_2 \cong 2\S_2\oplus3\S_3\oplus \S_{\bcirc}\oplus \S_{\posetV}\oplus \S_{\posetVop}\oplus U\mpoint$$
We note that both $\S_1$ and $\S_2$ are projective functors, hence so is $\S_1\otimes \S_2$ by Proposition~\ref{projective-tensor}.
}
\fresult

\bigskip\bigskip
\noindent\textbf{Acknowledgments.} The first author is grateful to Antoine Touz\'e for an inspiring talk about exponential functors, given at the Amiens group theory seminar, and to Robert Boltje, for fruitful discussions around Remark~\ref{uniqueness}.

\bigskip
\noindent
Serge Bouc, CNRS-LAMFA, Universit\'e de Picardie - Jules Verne,\\
33, rue St Leu, F-80039 Amiens Cedex~1, France.\\
{\tt serge.bouc@u-picardie.fr}

\medskip
\noindent
Jacques Th\'evenaz, Section de math\'ematiques, EPFL, \\
Station~8, CH-1015 Lausanne, Switzerland.\\
{\tt Jacques.Thevenaz@epfl.ch}

\end{document}